\documentclass[12pt]{amsart}
\usepackage[top=1in,bottom=1in,left=1in,right=1in]{geometry}
\usepackage{indentfirst}
\usepackage{amsfonts,amsmath,amsthm,amssymb}
\usepackage{longtable,enumerate,xcolor,seqsplit}
\usepackage{hyperref}

\newtheorem{theorem}{Theorem}[section]

\newtheorem{proposition}[theorem]{Proposition} 
\newtheorem{lemma}[theorem]{Lemma}

\newtheorem{remark}[theorem]{Remark}
\newtheorem{question}{Question}
\newtheorem{problem}[question]{Problem}
\newtheorem{hypothesis}[theorem]{Hypothesis}

\renewcommand{\leq}{\leqslant}
\renewcommand{\geq}{\geqslant}

\newcommand{\Soc}{{\rm Soc}}

\newcommand{\B}{\mathcal{B}}
\def\P{\mathcal{P}}

\newcommand{\Aut}{{\rm Aut}}

\newcommand{\PSL}{{\rm PSL}}
\newcommand{\PGL}{{\rm PGL}}

\newcommand{\PGGL}{{\rm P\Gamma L}}

\newcommand{\AGL}{{\rm AGL}}

\newcommand{\PSU}{{\rm PSU}}
\newcommand{\PGGU}{{\rm P\Gamma U}}
\newcommand{\PSp}{{\rm PSp}}

\newcommand{\la}{\lambda}

\DeclareMathOperator{\D}{\mathcal{D}}
\DeclareMathOperator{\C}{\mathcal{C}}

\begin{document}

\title{Analysing flag-transitive point-imprimitive $2$-designs}

\author{Alice Devillers  and Cheryl E. Praeger}\address{Address: Centre for the Mathematics of Symmetry and Computation, University of Western Australia, Perth, WA 6009, Australia.}
\email{alice.devillers@uwa.edu.au, cheryl.praeger@uwa.edu.au}
\thanks{This work was supported by an Australian Research Council Discovery Grant Project DP200100080.The authors thanks members of the Centre for the Mathematics of Symmetry and Computation at the University of Western Australia for their participation in early discussions of the project at the 2021 Research Retreat.}
\maketitle
\begin{abstract}
    In this paper 
we develop several general methods for analysing flag-transitive point-imprimitive $2$-designs, which give restrictions on both  the automorphisms and parameters of such designs. These  constitute a tool-kit for analysing these designs and their groups. We apply these methods to  complete the classification of flag-transitive, point-imprimitive $2$-$(v,k,\lambda)$ designs with $\lambda$ at most $4$.
\end{abstract}

\section{Introduction}
A \emph{$2$-$(v,k,\lambda)$ design} $\D=(\P,\B)$ consists of a set $\P$ of $v$ \emph{points} and a set $\B$ of \emph{blocks}
such that each block is a $k$-subset of $\P$ and each pair of distinct points is contained in $\lambda$ blocks.
To avoid degenerate cases we assume that $2<k<v$;  such designs are called \emph{nontrivial}. In general, the number of blocks $b:=|\B|$ is at least $v$ by Fisher's inequality (see \cite[1.3.8]{Dem}) and $\D$ is said to be \emph{symmetric} if $b=v$. We study (not necessarily symmetric) 2-designs $\D$ possessing a high degree of symmetry, namely they admit a subgroup $G$ of automorphisms  (permutations of $\P$ preserving $\B$) that acts transitively on the set of \emph{flags} (incident point-block pairs), and moreover leaves invariant a nontrivial partition of the point set $\P$, that is to say, $G$ is flag-transitive and point-imprimitive. We develop several general purpose methods for analysing flag-transitive point-imprimitive 2-designs, and then we test their effectiveness by applying them to  complete the classification of the flag-transitive, point-imprimitive $2$-$(v,k,\lambda)$ designs with $\lambda\leq 4$.

For a flag-transitive, point-imprimitive $2$-$(v,k,\lambda)$ design, the parameter $\lambda\ne 1$  by the celebrated work of Higman and McLaughlin~\cite[Proposition 3]{HM}. Moreover the classification of such designs with $\lambda=2$ was completed in \cite[Theorem 1.1]{DLPX}, showing that there are exactly two examples up to isomorphism. The classification for $\lambda\in\{3,4\}$ with $v<100$ was given in \cite[Theorem 2]{DP}, identifying nine designs up to isomorphism, and  
here we complete that work (solving \cite[Problem 3]{DP}) and proving that there are no examples with $100$ points or more. 

\begin{theorem}\label{t:lam34}
Let $\D=(\P,\B)$ be a  $2$-$(v,k,\lambda)$ design with $\lambda\leq 4$, which admits a flag-transitive, point-imprimitive subgroup of automorphisms. Then $v\in\{15,16,36,45,96\},$ and $\mathcal{D}$ is one of the eleven designs listed in \cite[Theorem 2]{DP}.
\end{theorem}

A number of papers have appeared recently providing classifications of flag-transitive, point-imprimitive $2$-designs under various parameter constraints. We give a brief survey of such results in Subsection~\ref{s:survey}, where we also introduce the relevant parameters we will use to describe the point-imprimitivity system, notably Hypothesis~\ref{h}.

\subsection{Survey and parameters}\label{s:survey}

As we mentioned above, if $G\leq \Aut(\D)$ is flag-transitive and point-imprimitive on a nontrivial $2-(v,k,\lambda)$ design $\D=(\P,\B)$, then $\lambda\geq 2$ by \cite[Proposition 3]{HM}. This condition was sharpened by Dembowski~\cite[2.3.7(a)]{Dem} in his 1968 book, where he proved that each of the following conditions must hold, where $r$ is the (constant) number of blocks containing a given point, and $(s,t)$ denotes the greatest common divisor $\gcd(s,t)$ for integers $s,t$.
\begin{equation}\label{e:dem}
(\lambda,r)\geq 2,\quad \lambda\leq (\lambda,r)((\lambda,r)-1),\quad (r-\lambda,k)\geq2,\quad  r\leq \lambda(k-3),\quad (v-1,k-1)\geq3,     
\end{equation}
We mention in passing that, very recently, the fourth inequality was strengthened by  Zhao and Zhou~\cite[Lemma 1.3]{ZZ} to $r\leq (r,\lambda)(k-3)$. After Dembowski's work the next significant breakthrough was due to Davies~\cite{Dav} in 1987. Davies showed by example that there are flag-transitive point-imprimitive designs with arbitrarily large $\lambda$, and also showed that, for a given $\lambda$, both the block-size $k$ and the number $v$ of points are bounded in terms of $\lambda$, and hence there are only finitely many flag-transitive point-imprimitive designs for each $\lambda$. Unfortunately Davies did not give upper bonds for $k, v$ as explicit functions of $\lambda$. 

The first such explicit bounds were due to O'Reilly-Regueiro~\cite{OR} in 2005 in the case of symmetric designs, for example, she showed that $k\leq \lambda(\lambda+1)$. These bounds were improved by Zhou and Praeger in \cite[Theorem 1.1]{PZ}, and were expressed in terms of additional parameters, namely the number $d$ of classes of a nontrivial point-partition $\mathcal{C}$,  the size $c$ of each class $\Delta\in\mathcal{C}$, and the (constant) size $\ell$ of each non-empty intersection $\Delta\cap B$ of a block $B$ and a class $\Delta$. The bounds were sufficiently good to show that for $\lambda\leq 10$, there are exactly $22$ feasible parameter tuples $(v,k,\lambda,c,d,\ell)$ meeting these bounds, \cite[Corollary 1.3 and Table 1]{PZ}. Very recently, in~\cite[Proposition 12]{MS}, Mandi\'{c} and \v{S}uba\v{s}i\'{c} made small improvements in the bounds for symmetric designs in \cite{PZ} (see also Proposition~\ref{p:MS}) and, building on classifications in \cite{LPR,P,PZ} of the examples with $\lambda\leq 4$, they were able to identify all examples with parameter tuples in  \cite[Table 1]{PZ} except for lines $10$ and $17$.

\begin{problem}\label{prob1}
Decide if there exist flag-transitive, point-imprimitive symmetric $2$-designs with parameter tuple $(v,k,\lambda,c,d,\ell)$ equal to $(288,42,6,8,36,2)$ or $(891,90,9,81,11,9)$, and if so classify all such designs.
\end{problem}

\noindent 
A resolution of Problem~\ref{prob1} would complete the classification of symmetric flag-transitive, point-imprimitive $2$-designs with $\lambda\leq 10$.
We note, in addition, that Montinaro~\cite{M} has recently classified all symmetric examples for which $k>\lambda(\lambda-3)/2$ and $\ell\geq3$.

For general (not necessarily symmetric) flag-transitive point-imprimitive designs with a given $\lambda$, explicit upper bounds for $k$ and $v$ were obtained in \cite[Theorem 1]{DP}. In this general case the upper bound on $k$ is cubic, namely $k\leq 2\lambda^2(\lambda-1)$;  the parametric restrictions obtained were sufficiently strong to list in  \cite[Proposition 8]{DP} all parameter sequences $(v,k,\lambda,c,d,\ell)$ meeting these restrictions with $\lambda\in\{3,4\}$, noting that the examples with $\lambda=2$ were classified in \cite[Theorem 1.1]{DLPX}.
In this paper, after discussing a number of general methods in Section~\ref{sec:methods}, we complete the classification of all examples arising from these parameter tuples in Section~\ref{ruleout}. 

We formalise the conditions we have been discussing in Hypothesis~\ref{h}. We will use these assumptions throughout the paper. A summary reference for the conditions in Hypothesis~\ref{h} may be found in  \cite[Section 2]{DP}. Note that we make no restrictions on the parameter $\lambda$ in Hypothesis~\ref{h}.

\begin{hypothesis}\label{h}  
\begin{enumerate}
    \item[(a)] $\mathcal{D}=(\mathcal{P}, \mathcal{B})$ is a $2$-$(v,k,\lambda)$ design, with point-set $\mathcal{P}$ of size $v$, block-set $\mathcal{B}$ of size $b$,  each block a $k$-subset of $\mathcal{P}$, and each point-pair lying in $\lambda$ blocks. Let $b=|\mathcal{B}|$, the number of blocks, and let $r$ be the number of blocks containing a given point.  
    
    \item[(b)] $\mathcal{D}$ admits a group $G$ of automorphisms such that $G$ is transitive on flags, and is imprimitive on points, preserving a point-partition $\mathcal{C}=\{\Delta_1,\dots, \Delta_d\}$ of size $d\geq2$ with classes $\Delta_i$ of size $c\geq 2$. Let $D=G^\mathcal{C}$ and $L=(G_{\Delta})^{\Delta}$ denote the induced permutation groups on $\mathcal{C}$ and $\Delta$, for $\Delta\in\mathcal{C}$. We may (see \cite[Theorem 5.5]{PS}), and will, assume that $G\leq L\wr D$. Let $K:=G_{(\mathcal{C})}$, the kernel of the $G$-action on $\mathcal{C}$.
    \item[(c)] For $B\in\mathcal{B}$ intersecting a class $\Delta\in\mathcal{C}$ nontrivally, the size  $\ell=|B\cap\Delta|$ is independent of the choices of $B, \Delta$; and $\ell\geq2$.
\end{enumerate}
\end{hypothesis} 

The general results presented in Section \ref{sec:methods} about flag-transitive, point-imprimitive  designs, include also various properties of the groups $L$ and $D$. Together these constitute a tool-kit of general methods for analysing these designs and their groups. 
We are in particular interested in the case $\lambda=3$ or $4$ with $v\geq100$ and  in Section \ref{ruleout} we apply these methods to prove Theorem~\ref{t:lam34}.

We present the methods in Section~\ref{sec:methods} to draw attention to the strategies we employ in the analysis in Section~\ref{ruleout}. We would be glad if this encouraged others to strengthen the methods in Section~\ref{sec:methods}, or obtain more extensive classifications of flag-transitive designs. In particular we wonder if these considerations might lead to improvements in the general upper bounds derived in \cite{DP}.

\section{General methods}\label{sec:methods}

In this section we assume that Hypothesis~\ref{h} holds and we use the notation introduced there. 
We begin by giving in Lemma~\ref{lem:parameters} a summary of basic restrictions on the parameters in Hypothesis~\ref{h}, which are in addition to the conditions in \eqref{e:dem}. Lemma~\ref{lem:parameters} is a simplification of the statements of  \cite[Lemmas 4 and 5]{DP}, together with an extra part (iii) which follows from an argument to be found in the proof of \cite[Proposition 9]{MS} (see also  Proposition \ref{p:MS} below). Our proof of part (iii) is essentially the first part of the proof of \cite[Proposition 9]{MS}.



\begin{lemma}\label{lem:parameters}
The following equalities, inequalities and divisibility conditions hold:
\begin{enumerate}[(i)]
\item $bk=vr$ and $r(k-1)=\lambda(v-1)$
\item $\ell\mid k\quad \mbox{and}\ 1<\ell < k$;
\item $\ell^2\mid c^2\lambda$
\item $\ell-1\leq \left(k-1-d(\ell-1)\right) \left(\ell-1\right)\leq \la-1$;
\item $\la(c-1)=r(\ell-1)$;
\item $k\mid \la\ell(\ell-1)^2d(d-1)$. 
\end{enumerate} 
\end{lemma}

\begin{proof}
Part (i)  is standard, see for example \cite[Lemma 4]{DP}. Parts (ii) and (v) are \cite[Lemma 5(ii) and (iii)]{DP}, respectively. For  part (iv), the first inequality folllows from the assertion in \cite[Lemma 5]{DP} that the quantity $x=k-1-d(\ell-1)$ is a positive integer, and the second inequality follows from   \cite[Lemma 5(viii)]{DP} on substituting for $x$. Part (vi) follows from   \cite[Lemma 5(vii)]{DP} on substituting for $x$.

Finally for part (iii), we choose distinct classes $\Delta,\Delta'\in\mathcal{C}$, and determine the cardinality of the set 
$\Pi:=\{(B,\alpha,\alpha')\mid B\in\mathcal{B}, \alpha\in B\cap\Delta, \alpha'\in B\cap\Delta'\}$. On the one hand there are $c^2$ choices for a pair $(\alpha,\alpha')\in\Delta\times\Delta'$ and each pair lies in exactly $\lambda$ blocks $B$, so $|\Pi|=c^2\lambda$. On the other hand each of the, say $n$, blocks $B$ meeting both $\Delta$ and $\Delta'$ nontrivially meets each of these classes in exactly $\ell$ points, so $|\Pi|=n\ell^2$. It follows that $c^2\lambda=n\ell^2$, and hence the number of these blocks $n=c^2\lambda/\ell^2$. Part (iii) follows since $n$ is an integer.
\end{proof}

\medskip
Next we derive  restrictions on the group $D$ in Hypothesis~\ref{h}.

\begin{proposition}\label{le:Druleout}
Assume that Hypothesis \ref{h} holds, and suppose that
$p$ is a prime such that:
\begin{enumerate}[(i)]
\item $0\leq d-p<\frac{k}{\ell}< p$, and
\item $p$ does not divide $b=\frac{vr}{k}$.
\end{enumerate}
Then $p$ does not divide $|D|$. 
\end{proposition}

\begin{proof}
Suppose  $p$ divides $|D|$ and $p$ satisfies the two stated conditions. Note the first condition implies $d<2p$.
Let $P$ be a Sylow $p$-subgroup of $G$.
Since $p$ divides $|D|$ and hence $|G|$, $P$ is not trivial. 
Moreover, the index of $G_B$ in $G$ is $b$, which is not divisible by $p$. Therefore there exists a block $B$ such that $P\leq G_B$.
Since $p$ divides $|D|$, $P$ acts non-trivially on $\mathcal{C}$. The condition $d<2p$ implies that $P^\mathcal{C}$ must have one orbit of size $p$ and $d-p$ fixed classes.
If $B$ contains a point in a class in the orbit of size $p$, then $k\geq \ell p$ since $P$ fixes $B$, contradicting $\frac{k}{\ell}< p$. 
Thus each point in $B$ lies in one of the $d-p$ fixed classes. This implies $k\leq \ell(d-p)$, contradicting $d-p<\frac{k}{\ell}$. 
This final contradiction proves that $p$ does not divide $|D|$.
\end{proof}

\medskip

\begin{proposition}\label{le:DBorbits} 
Assume that Hypothesis \ref{h} holds. Let $B\in\mathcal{B}$ and $X:=(G_B)^\mathcal{C}$, a subgroup of $D$. Also let $c_0$ be the (constant) length of the $K$-orbits in $\mathcal{P}$, and let $x:=c_0/\gcd(c_0,\ell)$.  Then $X$ has an orbit of length $k/\ell$ in $\mathcal{C}$ and
\begin{enumerate}
    \item[(a)] if $c_0=1$ (or equivalently, if $K=1$), then $|D:X|=b$; while
    
    \item[(b)] if $c_0>1$, then $c_0$ divides $c$, $x$ divides $b$, and $|D:X|$ divides $b/x$. 
\end{enumerate}
\end{proposition}

\begin{proof}
Since $G_B$ is transitive on $B$, $X$ is transitive on the set of $k/\ell$ classes intersecting $B$ non-trivially, proving the first assertion. Also, since $X=G_BK/K$ and $D=G/K$, it follows that $|D:X|=|G:G_BK|$. In particular, if $K=1$, or equivalently, $c_0=1$, then $|D:X|=|G:G_B|=b$ and part (a) holds. 

Assume now that $K\ne 1$, or equivalently that $c_0>1$. Since $K^\Delta\lhd\,G_\Delta^\Delta=L$, it follows that $c_0$ divides $c$. 
Since $G$ is block-transitive, its normal subgroup $K$ has orbits of equal length, say $b_0$,  in $\mathcal{B}$, so $b_0=|K:K_B|$ divides $b$. Let $\Delta\in\mathcal{C}$ be such that $B\cap\Delta\ne\emptyset$, and let $\alpha\in B\cap\Delta$. Since $G_B$ is transitive on $B$ it follows that $B\cap\Delta$ is a  $G_{B,\Delta}$-orbit of size $\ell$, and as $K_B\unlhd\, G_{B,\Delta}$, the $K_B$-orbits in $B\cap\Delta$ have equal length, say $\ell_0$, Thus $\ell_0 \mid  \ell$, and $|K:K_{B,\alpha}|=b_0\ell_0$ is divisible by $|K:K_\alpha|=c_0$.    Let $\ell':=\gcd(c_0,\ell_0)$. Then $\ell'\mid \gcd(c_0,\ell)$ and $c_0/\ell'$ divides $b_0$. Since 
\[
\frac{c_0}{\ell'} = \frac{c_0}{\gcd(c_0,\ell)}\cdot \frac{\gcd(c_0,\ell)}{\ell'} = x\cdot \frac{\gcd(c_0,\ell)}{\ell'}
\]
it follows that $x$ divides $b_0$ which, in turn, divides $b$. Finally 
\[
|D:X|=|G:G_BK|=\frac{|G:G_B|}{|G_BK:G_B|}= \frac{|G:G_B|}{|K:K_B|}= \frac{b}{b_0}= \frac{b}{x}\cdot \frac{1}{b_0/x}, 
\]
and it follows that $|D:X|$ divides $b/x$, and part (b) holds, completing the proof. 
\end{proof}

\medskip
Now we turn our attention to  restrictions on the group $L$ in Hypothesis~\ref{h}.

\begin{proposition}\label{le:L} Assume that Hypothesis \ref{h} holds, and let $\Delta\in\mathcal{C}$ and $\alpha,\beta\in\Delta$ be distinct points.
\begin{enumerate}
    \item[(a)]  If  $\ell=2$, then $L$ is $2$-transitive on $\Delta$ of degree $c$.
    \item[(b)] If $\ell\geq3$, then $L_{\{\alpha,\beta\}}$ has an orbit in $\Delta\setminus\{\alpha,\beta\}$ of size at most $(\ell-2)\lambda$. In particular, if $c>2+(\ell-2)\lambda$ then $L$ is not $3$-transitive.
\end{enumerate}
\end{proposition}

\begin{proof}
(a) Assume that $\ell=2$ and let $(\alpha',\beta')$ be a second pair of distinct points of $\Delta$. Let $B, B'\in\mathcal{B}$  be such that $\alpha,\beta\in B$ and $\alpha', \beta'\in B'$. Since $G$ is flag-transitive, there exists $g\in G$ mapping $(\alpha,B)$ to $(\alpha',B')$. Then since $\alpha'=\alpha^g\in \Delta\cap \Delta^g$ it follows that $g$ fixes $\Delta$, so $g^\Delta\in L$. Also, since $B^g=B'$ and $\ell=2$, the element $g$ must map $\beta$ to $\beta'$. It follows that $L$ is 2-transitive.

(b) Let $\mathcal{B}(\alpha,\beta)=\{B_1,\dots, B_\lambda\}$ be the set of $\lambda$ blocks
containing $\{\alpha,\beta\}$. Then  $G_{\{\alpha,\beta\}}$ leaves $\mathcal{B}(\alpha,\beta)$ 
invariant, and hence leaves invariant the union $X:=\cup_{i=1}^\lambda B_i\cap \Delta$. Thus  
$X\setminus \{\alpha,\beta\}$ is a subset of $\Delta\setminus\{\alpha,\beta\}$ of size at most $(\ell-2)\lambda$, it is non-empty since $\ell\geq 3$, and it is preserved by $G_{\{\alpha,\beta\}}$. Any 
$G_{\{\alpha,\beta\}}$-orbit in $X\setminus \{\alpha,\beta\}$  (which is also an $L_{\{\alpha,\beta\}}$-orbit) has size at most $(\ell-2)\lambda$. 
If $L$ were 3-transitive, then $L_{\{\alpha,\beta\}}$ would be transitive on $\Delta\setminus\{\alpha,\beta\}$, and hence $c\leq 2+(\ell-2)\lambda$.
\end{proof}

\medskip
We now explore properties of certain $p$-subgroups of $G$.
\begin{proposition}\label{le:plam}
Assume that Hypothesis \ref{h} holds, and suppose that $p$ is a prime dividing $|G|$ such that $p$ does not divide $\lambda$. Suppose also that a nontrivial $p$-subgroup $P$ of $G$ fixes at least two distinct points $\alpha, \beta$, and let $\Delta\in\mathcal{C}$ such that $\alpha\in\Delta$.
\begin{enumerate}
    \item[(a)] Then $P$ fixes at least one block $B$ containing $\{\alpha,\beta\}$ and hence leaves invariant the $\ell$-subset $B\cap\Delta$ of $\Delta$ and the $k/\ell$-subset $\mathcal{C}(B) :=\{\Delta'\in\mathcal{C}\mid B\cap \Delta'\ne\emptyset\}$ of $\mathcal{C}$. 
    \item[(b)] Moreover, if $p>\lambda$ then $P$ fixes setwise each block containing $\{\alpha,\beta\}$.
\end{enumerate}
\end{proposition}

\begin{proof}
Let $\mathcal{B}(\alpha,\beta)$ be the set of $\lambda$ blocks
containing $\{\alpha,\beta\}$. Then  $P$ leaves $\mathcal{B}(\alpha,\beta)$ invariant. Since $P$ is a $p$-group, all its orbits on $\mathcal{B}(\alpha,\beta)$ are $p$-powers. Since $p$ does not divide $\lambda$, it follows that $P$ fixes at least one block $B\in\mathcal{B}(\alpha,\beta)$. Part (a) follows.
It also follows that if $p>\lambda$ then $P$ must fix setwise each block in $\mathcal{B}\{\alpha,\beta\}$.
\end{proof}




\medskip

Our next results concern a \emph{subdesign} of a design $\mathcal{D}=(\mathcal{P}, \mathcal{B})$. This is a pair $\mathcal{D}_0=(\mathcal{P}_0, \mathcal{B}_0)$ where $\mathcal{P}_0\subset \mathcal{P}$ and $\mathcal{B}_0$ is a collection of $k_0$-subsets of $\mathcal{P}_0$ such that each $B_0\in\mathcal{B}_0$ is contained in some block of $ \mathcal{B}$; and such that $\mathcal{D}_0$ is a $2$-design (or sometimes only a $1$-design). A subdesign is called \emph{complete} if $k_0=k$, and in this case $\mathcal{B}_0\subset \mathcal{B}$, see for example \cite[p.\,31]{HP}. Subdesigns provide useful information about the structure of a design. For example the design of points and lines of a projective space $\mathrm{PG}_n(q)$ is a flag-transitive $2-(\frac{q^n-1}{q-1},q+1,1)$ design and for each proper subspace, the set of lines it contains forms a complete subdesign. Similar complete subdesigns arise for the designs of points and lines of a Desarguesian affine space, and in both cases these subdesigns inherit a flag-transitive action from the automorphism group of the original design. In these examples the designs are point-primitive. 

\begin{question}
Do there exist $2$-designs admitting a flag-transitive, point-imprimitive group of automorphisms and containing a complete subdesign?
\end{question}

In Proposition~\ref{cor:subdesign} we show that under additional assumptions, there is a point-transitive complete subdesign associated with a Sylow subgroup of $K$. This gives a strong restriction on the parameters which we exploit in Section~\ref{ruleout} to exclude certain parameter tuples. However we do not know of any examples satisfying Hypothesis \ref{h} where conditions (a)--(c) of Proposition~\ref{cor:subdesign} hold.

\begin{proposition}\label{cor:subdesign} Assume that Hypothesis \ref{h} holds and that there exists a prime $p$ satisfying the following conditions:
\[
\mbox{(a)\ \ $p$ divides $|K|$,\quad (b)\ \ $p$ does not divide $c$, \quad (c)\ \ $p>\lambda$.}
\]
Then $\mathcal{D}$ has a complete point-transitive subdesign which is  a $2-(df,k,\lambda)$ design, where $f$ satisfies 
\[
\mbox{$\ell< f<c$\quad and\quad $c\equiv f\pmod p$.}
\]
In particular, $k-1$ divides $\lambda(df-1)$.
\end{proposition}

\begin{proof}
Let $P$ be a Sylow $p$-subgroup of $K$, and let $\mathcal{P}_0={\rm fix}_\mathcal{P}(P)$ and $\mathcal{B}_0 =\{ B\in\mathcal{B}\mid B\subseteq \mathcal{P}_0\}$.
We show that  $(\mathcal{P}_0, \mathcal{B}_0)$ is a complete subdesign of $\mathcal{D}$, and is a  $2-(df,k,\lambda)$ design for some integer $f$ as in the statement, and that $N_G(P)$ induces a point-transitive automorphism group. 

By condition (a), $P$ is non-trivial. We see that $G=N_G(P)K$ by a `Frattini argument': for $g\in G$, $P^g$ is a Sylow $p$-subgroup of $K^g=K$ and hence $P^g=P^x$ for some $x\in K$, whence $y:=gx^{-1}\in N_G(P)$ so $g=yx\in N_G(P)K$. It follows that $N_G(P)^\mathcal{C} =G^\mathcal{C}=D$, and in particular $N_G(P)$ is transitive on $\mathcal{C}$. 

Now $P^{\Delta}\ne 1$ for some $\Delta\in\C$, since $P\ne 1$. 
By condition (b), the set ${\rm fix}_\Delta(P)$ of fixed points of $P$ in $\Delta$ is non-empty, and is a proper subset of $\Delta$ since $P^{\Delta}\ne 1$. Hence $f:=|{\rm fix}_\Delta(P)|$ satisfies $0<f<c$  and $c\equiv f\pmod p$. Let $\alpha\in {\rm fix}_\Delta(P)$ and note that $P$ is a Sylow $p$-subgroup of $K_\alpha$. If $g\in G$ is such that $P^g\leq G_\alpha$, then $P^g$ is contained in $K^g\cap G_\alpha=K_\alpha$ so $P^g$ is also a Sylow $p$-subgroup of $K_\alpha$ and hence $P^g=P^x$ for some $x\in K_\alpha$. Thus $P^g$ is conjugate to $P$ in $G_\alpha$ and hence,  by \cite[Corollary 2.24]{PS}, $N_G(P)$ is transitive on ${\rm fix}_\mathcal{P}(P)$. It follows that $P$ fixes the same number $f$ of points in each class of $\mathcal{C}$, and hence $|\mathcal{P}_0|=|{\rm fix}_\mathcal{P}(P)|=df$.

Let $\alpha,\beta\in {\rm fix}_\mathcal{P}(P)$ with $\alpha\ne\beta$, and let $B$ be a block containing $\{\alpha,\beta\}$.
By Proposition \ref{le:plam}(b), $P$ fixes $B$ setwise. Hence $P$ fixes setwise each non-trivial intersection $B\cap\Delta'$, for $\Delta'\in\mathcal{C}$, and each such intersection has size $\ell$. Now  $\lambda\geq \ell$ by Lemma \ref{lem:parameters}(iv)
, and therefore, by condition (c), $p>\ell$. Thus each non-trivial intersection $B\cap\Delta'$ is fixed pointwise by $P$, and hence $f\geq \ell$ and $P$ fixes $B$ pointwise, that is, $B\subseteq {\rm fix}_\mathcal{P}(P)$. This implies that each of the $\lambda$ blocks containing $\{\alpha,\beta\}$ lies in $\mathcal{B}_0$, and it follows that  $(\mathcal{P}_0, \mathcal{B}_0)$ is a $2-(df,k,\lambda)$ design, and hence a complete subdesign of $\mathcal{D}$, admitting $N_G(P)$ acting as a point-transitive automorphism group. In particular  $k-1$ divides $\lambda(df-1)$, see for example \cite[Lemma 4(i)]{DP}.

It remains to prove that $f>\ell$. Suppose to the contrary that $f=\ell$. Then whenever $B\cap\Delta\ne\emptyset$ we have $B\cap\Delta = {\rm fix}_{\Delta}(P)$ of size $\ell$. Thus each block $B$ containing at least two points of ${\rm fix}_\mathcal{P}(P)$ is the disjoint union of the subsets ${\rm fix}_\Delta(P)$ over all classes $\Delta\in\mathcal{C}$ such that $B\cap\Delta\ne\emptyset$. 
Fix $\Delta\in\mathcal{C}$, choose distinct points $\alpha,\beta\in {\rm fix}_\Delta(P)$, and let $N$ be the number of pairs $(B,\Delta')$ such that $B\in\mathcal{B}$, $\Delta'\in\mathcal{C}\setminus\{\Delta\}$, and $B\cap \Delta$, $B\cap \Delta'$ are both non-empty.  The blocks $B$ occurring are precisely the $\lambda$ blocks containing $\{\alpha,\beta\}$ and each such block $B$ occurs in pairs $(B,\Delta')$ for exactly $k/\ell-1$ classes $\Delta'$, so $N=\lambda(k/\ell-1)$. On the other hand, for each of the $d-1$ classes $\Delta'\ne\Delta$, choosing $\gamma\in {\rm fix}_{\Delta'}(P)$, we see that $\Delta'$ occurs in pairs $(B,\Delta')$ for precisely the $\lambda$ blocks $B$ containing $\{\alpha,\gamma\}$, so $N=\lambda(d-1)$. We conclude that $k/\ell=d$, However, this implies that each of the $\lambda$ blocks containing $\{\alpha,\beta\}$ is equal to ${\rm fix}_{\mathcal{P}}(P)$ which, in turn, implies that $\lambda=1$, contradicting \cite{HM}. Thus $f>\ell$.
\end{proof}

\medskip
 There are also flag-transitive subdesigns (not complete ones) and `quotient designs' arising from flag-transitive, point-imprimitive designs, see Proposition~\ref{p:MS} below and the design construction in \cite[Construction 3.1]{CP}. Their natural definition involves the possibility of `repeated blocks', that is distinct blocks incident with exactly the same subset of points, so we use a more formal incidence structure in their definition. Given the assumptions in Hypothesis~\ref{h}, let $\Delta\in\mathcal{C}$, and define $\D(\Delta)=(\Delta,\B(\Delta),\mathcal{I}(\Delta))$, where 
 \[
 \B(\Delta)=\{ B\cap\Delta\mid B\in\B, B\cap\Delta\ne \emptyset\},\ \mbox{and}\ 
 \mathcal{I}(\Delta) = \{(\alpha,B\cap\Delta)\mid \alpha\in B\cap\Delta, B\cap\Delta\in B(\Delta)\},
 \]
and $\D(\mathcal{C})=(\mathcal{C},\mathcal{C}(\mathcal{\B}),\mathcal{I}(\mathcal{C}))$, with 
$\mathcal{C}(\mathcal{\B})$ the set of all $\mathcal{C}(B)$ for $B\in\B$, where
 \[
\mathcal{C}(B)=\{ \Delta\in\mathcal{C}\mid B\cap\Delta\ne \emptyset\},\ \mbox{and}\ 
 \mathcal{I}(\mathcal{C}) = \{(\Delta,\mathcal{C}(B))\mid \Delta\in \mathcal{C}(B)\}.
 \]
The notion of a $2$-design carries over to these point-block incidence structures (each pair of points incident with the same number of blocks), and automorphisms are permutations of the point set and the block set which preserve the incidence relation. The group $G_\Delta$ or $G$  acts as a flag-transitive group of automorphisms on $\D(\Delta)$ or  $\D(\mathcal{C})$, respectively. This means in particular that, for each design, each block occurs with the same \emph{block-multiplicity} (the number of blocks incident with the same subset of points). In the following proposition, part (a) was proved in both  \cite[Propositions 6 and 8]{MS} and \cite[Corollary 2.2 and Theorem 2.3]{M}; and part (b) is proved in   \cite[Propositions 7 and 9]{MS} (and in both parts the block multiplicities may be greater than $1$). 
In both papers these observations were applied to strengthen the results of \cite{PZ}
(see  \cite[Propositions 12 and 9]{MS} and \cite[Theorem 2.4]{M}).
Note that although the applications in \cite{MS,M} are to symmetric designs, Proposition~\ref{p:MS} is valid for all flag-transitive point-imprimitive designs.

\begin{proposition}\label{p:MS}
Assume that Hypothesis \ref{h} holds, let $\Delta\in\mathcal{C}$, and consider the incidence structures $\D(\Delta)=(\Delta,\B(\Delta),\mathcal{I}(\Delta))$ and $\D(\mathcal{C})=(\mathcal{C},\mathcal{C}(\mathcal{\B}),\mathcal{I}(\mathcal{C}))$ as defined above. Then 
\begin{enumerate}
    \item[(a)]  $\D(\Delta)$ is a $2-(c,\ell,\lambda)$ design with $L=G_\Delta^\Delta$ acting flag-transitively, and the block multiplicity $\theta$  equals 
    $\theta=|G_{\alpha,B\cap\Delta}:G_{\alpha,B}|$ and divides $\lambda$, where $B\cap\Delta\in\B(\Delta)$ and $\alpha\in B\cap\Delta$; 
    
    \item[(b)] $\D(\mathcal{C})$ is a $2-(d, k/\ell,c^2\lambda/\ell^2)$ design with $D=G^\mathcal{C}$ acting flag-transitively.  
\end{enumerate}
\end{proposition}

\section{Ruling out designs with $\lambda$ small}\label{ruleout}

In this section we assume that Hypothesis~\ref{h} holds with $\lambda=3$ or $4$ and $v\geq100$. 
In \cite[Proposition 8]{DP}, it was shown that the parameter tuple $(\lambda, v, k, r, b, c, d, \ell)$ must belong to a list of 18 possibilities.  However two of those tuples should not have been listed since the value of $vk/r=b$ for them was not an integer  (namely the lines in \cite[Proposition 8]{DP} with $(\lambda, v)=(4, 196)$ and $(4, 435)$ should have been omitted). Moreover one of the remaining tuples does not satisfy Lemma \ref{lem:parameters}(iii), namely the tuple $(\lambda, v,c,d,\ell)=(3,561,17,33,2)$. 
Thus we have $15$ parameter tuples to consider, and we list these in Table~\ref{t:remainingcases}.  We will apply the results from Section \ref{sec:methods} to show that there are no designs for any of these parameter tuples. In the last column of Table~\ref{t:remainingcases} we record  the results in this section which rule out each of these parameter tuples.
In our statements, Line numbers refer to Table~\ref{t:remainingcases}.
Moreover, in our proofs we frequently mention use of the computational algebra package Magma~\cite{magma}, and we comment briefly in Remark~\ref{rem:magma}  on the main ways in which this is used.


\begin{table}[]
    \centering
\begin{tabular}{|c|cccccccc|l|}
\hline
Line&$\la$&$v$&$k$&$r$&$b$&$c$&$d$&$\ell$&Result\\
\hline
1&$3$&$100$&$12$&$27$&$3^2\cdot 5^2$&$10$&$10$&$2$& Lemma~\ref{cor:lines1and8and11}\\
2&$3$&$120$&$18$&$21$&$2^2\cdot 5\cdot 7$&$8$&$15$&$2$& Lemma~\ref{cor:plam}\\
3&$3$&$120$&$18$&$21$&$2^2\cdot 5\cdot 7$&$15$&$8$&$3$& Lemma~\ref{cor:plam}\\
4&$3$&$256$&$18$&$45$&$2^7\cdot 5$&$16$&$16$&$2$&Lemma~\ref{cor:DB}\\
5&$3$&$561$&$36$&$48$&$2^2\cdot 11\cdot 17$&$33$&$17$&$3$& Lemma~ \ref{cor:conditionL}\\
6&$3$&$1156$&$36$&$99$&$11\cdot 17^2$&$34$&$34$&$2$& Lemma~ \ref{cor:easy}\\
\hline
7&$4$&$100$&$12$&$36$&$2^2\cdot 3\cdot 5^2$&$10$&$10$&$2$& Lemma~\ref{cor:lines1and8and11}\\
8&$4$&$231$&$24$&$40$&$ 5\cdot 7\cdot 11$&$11$&$21$&$2$& Lemma~\ref{cor:line9}\\
9&$4$&$231$&$24$&$40$&$ 5\cdot 7\cdot 11$&$21$&$11$&$3$& Lemma~ \ref{cor:DB}\\
10&$4$&$280$&$32$&$36$&$3^2\cdot 5\cdot 7$&$10$&$28$&$2$& Lemma~\ref{cor:lines1and8and11}\\
11&$4$&$280$&$32$&$36$&$3^2\cdot 5\cdot 7$&$28$&$10$&$4$& Lemma~\ref{cor:line12}\\
12&$4$&$484$&$24$&$84$&$2\cdot 7\cdot 11^2$&$22$&$22$&$2$& Lemma~ \ref{cor:DB} \\
13&$4$&$1976$&$80$&$100$&$2\cdot 5\cdot 13\cdot 19$&$26$&$76$&$2$& Lemma~ \ref{cor:easy}\\
14&$4$&$1976$&$80$&$100$&$2\cdot 5\cdot 13\cdot 19$&$76$&$26$&$4$& Lemma~ \ref{cor:conditionL}\\
15&$4$&$2116$&$48$&$180$&$3\cdot 5\cdot 23^2 $&$46$&$46$&$2$& Lemma~ \ref{cor:easy}\\
\hline
\end{tabular}
    \caption{Remaining parameter sets for flag-transitive imprimitive designs with $\lambda\leq 4$}
    \label{t:remainingcases}
\end{table}    

 \begin{remark}\label{rem:magma}
 {\rm 
We make some comments here on our use of computation in analysing the possibilities for flag-transitive $2-(v,k,\lambda)$ designs, for some specified parameters  $v, k,\lambda$, noting  that these parameters determine the number $b$ of blocks since $bk(k-1)=v(v-1)\lambda$. 
 
\medskip\noindent
(a)  At many stages in our proofs we have information about a primitive permutation group of a given degree (usually $c$ or $d$ where $v=cd$), and we use the library of primitive groups in Magma~\cite{MAGMA1} to determine a list of possibilities for the group with the required properties. 

\medskip\noindent (b)  Occasionally in our proofs we have specified the point set $\mathcal{P}$ of size $v$, and a transitive  permutation group $G\leq{\rm Sym(}\mathcal{P})$ as a candidate for the flag-transitive group $G$. Moreover we have specified a candidate subgroup $H$ of  index $b$ in $G$ for the stabiliser $H=G_B$ of a block $B$. Since $G$ should be flag-transitive, the block $B$ should be an $H$-orbit in $\mathcal{P}$ of length $k$.  We use Magma~\cite{magma} to identify all $H$-orbits of length $k$ in $\mathcal{P}$. These will be the possibilities for the block $B$ left invariant by $H$. 
 If the set $\mathcal{B}$ of $G$-images of such a subset $B$ forms the block set of a $2$-design, then the group $G$ will act flag-transitively on this design, giving an example. To determine whether this is the case for a given $H$-orbit $B$ of length $k$, we need to confirm whether or not each pair of points is contained in exactly $\lambda$ blocks in $\mathcal{B}$. Since $G$ is transitive on $\mathcal{P}$ it is sufficient to test this for point-pairs $\{\alpha,\beta\}$, where $\alpha\in\mathcal{P}$ is a fixed chosen point, and $\beta$ runs over a set of representative  points, one from each $G_\alpha$-orbit in $\mathcal{P}\setminus\{\alpha\}$.This check is carried out computationally 
using~\cite{magma}.
}
 \end{remark}

First we observe that in each Line the groups $L, D$ are primitive permutation groups.

\begin{lemma}\label{le:prim}
If Hypothesis \ref{h} is satisfied with $\lambda\in\{3,4\}$, then $D$ and $L$ are primitive.
\end{lemma}

\begin{proof}
Note that, in Table~\ref{t:remainingcases}, there is no pair of
parameter tuples with the same $(\la,v,k)$ and with the class size $c$ in one tuple a multiple of the class size in the other. It follows that $L$ is primitive of degree $c$, and also $D$ is primitive of degree $d$.  
\end{proof}

\medskip
Proposition \ref{le:Druleout} directly allows us to rule out three Lines  of Table~\ref{t:remainingcases}.

\begin{lemma}\label{cor:easy}
There are no examples with parameter tuple in Line $6, 13$ or $15$ of Table~\ref{t:remainingcases}.
\end{lemma}

\begin{proof} 
Assume that Hypothesis~\ref{h} holds with parameter tuple as in one of Lines 6, 13 or 15 of Table~\ref{t:remainingcases}. 
By Lemma~\ref{le:prim}, $D$ is primitive of degree $d=34,76,46$ respectively. Checking the possibilities for such groups using Magma~\cite{MAGMA1}, we find that  $D=A_d$ or $S_d$.
Let $p$ be the largest prime less than $d$.  
Then $p$ divides $|D|$ and also $p$ satisfies conditions (i) and (ii) of Proposition \ref{le:Druleout} (see details in Table \ref{t:easy}). This contradicts Proposition~\ref{le:Druleout}. Thus there are no examples.

\begin{table}[]
    \centering
\begin{tabular}{|c|cccc|}
\hline
Line&$d-p$&$k/\ell$& $p$&$b$\\
\hline
6&$3$&$18$&$31$&$ 11\cdot 17^2$\\
13&$3$&$40$&$73$&$2\cdot 5\cdot 13\cdot 19$\\
15&$3$&$24$&$43$&$3\cdot 5\cdot 23^2$\\
\hline
\end{tabular}
    \caption{Table for the proof of Lemma~\ref{cor:easy}}
    \label{t:easy}
\end{table}
\end{proof}

\medskip
Next we use Proposition \ref{le:DBorbits} to rule out three more Lines  of Table~\ref{t:remainingcases}.

\begin{lemma}\label{cor:DB}
There are no examples with parameter tuple in Line 4, 9 or 12 of Table~\ref{t:remainingcases}.
\end{lemma}
\begin{proof}
The argument for Lines 9 and 12 requires only the first assertion of Proposition~\ref{le:DBorbits}, so we deal with these Lines first. The parameters $d, k/\ell, b$ in Lines 9 and 12 are given in Table~\ref{t:DB}.
There are eight primitive groups of degree $11$ (Line 9), and four primitive groups of degree $22$ (Line 12). An exhaustive search (with Magma~\cite{MAGMA1}) of all subgroups of index dividing $b$ for each of these primitive groups shows that none has an orbit of size $k/\ell$.
This contradicts Proposition~\ref{le:DBorbits}.
\
\begin{table}[]
    \centering
\begin{tabular}{|c|ccc|}
\hline
Line&$d$&$k/\ell$& $b$\\
\hline
9 &$11$ &$8$  &$5\cdot 7\cdot 11$\\
12 &$22$ &$12$ &$2\cdot 7 \cdot 11^2$\\
\hline
\end{tabular}
    \caption{Table for the proof of Lemma~\ref{cor:DB}}
    \label{t:DB}
\end{table}

Now we consider  Line 4  of Table~\ref{t:remainingcases}, where we have $(c,d,\ell,k/\ell)=(16,16,2, 9)$, and $b=2^7\cdot 5=640$. In particular $5$ divides $|G|$ since $G$ is block-transitive. 
%
Note that $p=5> \lambda =3$, and that $p$ does not divide $c$. 
Suppose first that $5$ divides $|K|$, so that Proposition~\ref{cor:subdesign} applies with $p=5$.
Then there exists an integer $f$ such that $\ell< f<c$, $c\equiv f\pmod p$, and $k-1$ divides $\lambda(df-1)$, which is a contradiction.
Thus $|K|$ is not divisible by $5$.
Since $5$ divides $|G|=|D|\cdot|K|$, it follows that $5$ divides $|D|$.
 There are $18$ primitive groups $D$ of degree $d=16$ with order divisible by $5$, so $D$ must be one of these by Lemma \ref{le:prim}. Let $X=(G_B)^\mathcal{C}\leq D$. By Proposition \ref{le:DBorbits}, $X$ has an orbit of size $9$.

If $K\neq 1$ then, since $L$ is primitive, the $K$-orbits in $\mathcal{P}$ have size $c=16$, so $|D:X|$ divides $b/8=2^3\cdot 5=80$  by Proposition \ref{le:DBorbits}. 
However, an exhaustive search using Magma~\cite{MAGMA1} shows that none of the $18$ possibilities for $D$ has a subgroup of index dividing $80$ with an orbit of size $9$.
Thus $K=1$, so $G\cong D$ has order divisible by the number $f=bk=640\cdot 18$ of flags, and the index  $|D:X|=640$. 
Again a computer search using Magma~\cite{MAGMA1}  shows that the only  primitive group $D$ of degree $d=16$ such that $|D|$ is divisible by $f$, and $D$ has a subgroup of index $640$ with an orbit of size $9$, is the affine group $[2^4]:S_6$. Thus $G=[2^4]:S_6$ and $G_\Delta=S_6$. However $G_\alpha$ is a subgroup of $G_\Delta$ with index $c=16$, and $S_6$ has no such subgroup, which is a contradiction.
\end{proof}

\medskip

We rule out two more  Lines  of Table~\ref{t:remainingcases} with Proposition \ref{le:L}.

\begin{lemma}\label{cor:conditionL}
There are no examples with parameter tuple in Line 5 or 14 of Table~\ref{t:remainingcases}.
\end{lemma}
\begin{proof}
In Line 5, $\ell=3$, and $L$ is primitive of degree $c=33$ by Lemma \ref{le:prim}, but $L$ is not 3-transitive by Proposition \ref{le:L}, since $c=33>2+(\ell-2)\lambda=5$. However each primitive group of degree $c=33$, namely $\PSL(2,32),\PGGL(2,32), A_{33}$, or $S_{33}$, is 3-transitive. In Line 14, $\ell=4$, $c=76$, and the only primitive groups of degree $c$ are $A_{76}$ and $S_{76}$, so $L$ is 3-transitive, contradicting Proposition \ref{le:L}, since $c=76>2+(\ell-2)\lambda=10$. 
\end{proof}

\medskip
In the next lemma we rule out Lines  2 and 3  of Table~\ref{t:remainingcases} using a combination of Propositions~\ref{le:DBorbits} and~\ref{cor:subdesign}.

\begin{lemma}\label{cor:plam}
There are no examples with parameter tuple in Line 2 or 3 of Table~\ref{t:remainingcases}.
\end{lemma}
\begin{proof}
In Line 2 or 3 we have $(c,d,\ell,k/\ell)=(8,15,2, 9)$ or $(15,8,3,6)$, respectively, and $b=2^2\cdot5\cdot7=140$. In particular $7$ divides $|G|$. 
%
Note that $p=7> \lambda =3$, and that $p$ does not divide $c$. 
If $7$ divides $|K|$, then Proposition~\ref{cor:subdesign} applies with $p=7$, so
 there exists an integer $f$ such that $\ell< f<c$, $c\equiv f\pmod p$, and $k-1$ divides $\lambda(df-1)$. However there is no such integer in either case. 
Hence $|K|$ is not divisible by $7$.
Since $7$ divides $|G|=|D|\cdot|K|$, it follows that $7$ divides $|D|$.

First consider Line 2. By Lemma \ref{le:prim}, $D$ is a primitive group of degree $d=15$, and also $|D|$ is divisible by $7$. There are four such groups, namely $A_7$, $\PSL(4,2)$, $A_{15}$ and $S_{15}$, so $D$ is one of these. Let $X=(G_B)^\mathcal{C}\leq D$. By Proposition \ref{le:DBorbits}, $X$ has an orbit of size $9$.
If $K\neq 1$ then, since $L$ is primitive, the $K$-orbits in $\mathcal{P}$ have size $c=8$, so $|D:X|$ divides $140/4=35$  by Proposition \ref{le:DBorbits}. 
However, an exhaustive search using Magma~\cite{MAGMA1}  shows that none of the four possibilities for $D$ has a subgroup of  index dividing $35$ with an orbit of size $9$.
Thus $K=1$ and  $|D:X|=140$. Again a computer search shows that  the only primitive group of degree $d=15$ that has a subgroup of index $140$ with an orbit of size $9$ is $A_7$, so $D=A_7$.

Thus $G\cong D=A_7$, $G_\Delta=\PSL(2,7)$ and, for $\alpha\in\Delta$, $G_\alpha$ is a Frobenius subgroup $F_{21}$ of order $21$. Now $A_7$ has a unique conjugacy class of subgroups $F_{21}$, and we may identify the point-set $\mathcal{P}$  with the set of right cosets of such a subgroup $G_\alpha$. A block stabiliser $G_B$  is a subgroup of index 140 with an orbit of size $k=18$ on points. Using \cite{magma}, we checked that there is a unique conjugacy class of subgroups of $G=A_7$ of index 140 and that a subgroup $G_B$ in this class has six orbits of length $18$ on points. However, a further computation with \cite{magma} showed that, for each of these $18$-subsets $B$, the set of $G$-images of $B$ is not the block-set of a $2$-design.

Consider now Line 3. As we showed above, $|D|$ is divisible by $7$ while $|K|$ is not divisible by $7$. Thus a Sylow $7$-subgroup $P$ of $G$ is non-trivial and $P\cong P^\mathcal{C}$,  a Sylow $7$-subgroup of $D$. Since $d=8$, $P$ fixes  one class $\Delta\in\mathcal{C}$ and acts transitively on the other $7$ classes. 
Suppose that  $7$ does not divide $|L|$. Then $P$ fixes $\Delta$ pointwise and,  by Proposition~\ref{le:plam}(b), $P$ fixes setwise each block intersecting $\Delta$ nontrivially.  Let $B$ be such a block.
 Then $P$ fixes setwise the $6$-subset $\mathcal{C}(B)$ of $\mathcal{C}$ (defined in Proposition~\ref{le:plam}), a contradiction. 
Hence $7$ divides $|L|$.
 
 Thus $L$ is a primitive subgroup of $S_{15}$ with order divisible by $7$, and $L$ is not 3-transitive by Proposition \ref{le:L}(b), since $c=15>2+(\ell-2)\lambda=5$. Checking with \cite{MAGMA1}, we see that $L$ must be either $A_7$ or $\PSL(4,2)$, and in particular $L$ is simple. Since $K^\Delta\lhd L$ and $7$ does not divide $|K|$, it follows that $K=1$, and hence $G\cong D\leq S_8$.  Thus $8!\geq |G|=8\cdot|G_\Delta|\leq 8\cdot |L|$ and it follows that $L=A_7$ and $G\cong D=A_8$. By Proposition \ref{le:DBorbits}, $D$ has an index $140$ subgroup $X$ with an orbit of size $6$. This however implies that the setwise stabiliser in $D$ of a $6$-element subset (which is isomorphic to $S_6$ and has index $28$ in $D$) must have a subgroup of index $140/28=5$.  This contradiction completes the proof.
\end{proof}

 \medskip
There are just five lines of Table~\ref{t:remainingcases} still to be resolved, namely Lines 1, 7, 8, 10, 11. We consider these together to obtain restrictions on the groups $D$ and $K$. For a group $H$ and prime $p$, $O_p(H)$ denotes the largest normal $p$-subgroup of $H$.



\begin{table}
\begin{center}
\begin{tabular}{|c|llllll|p{10.5cm}|}
\hline
Line&$\lambda$&$k$&$b$&$c$&$d$&$\ell$&possible groups $D$\\
\hline
1&$3$&$12$&$3^2\cdot 5^2$&$10$&$10$&$2$&$A_5,S_5,\PSL(2,9),S_6$ (out of $9$)\\
7&$4$&$12$&$2^2\cdot 3\cdot 5^2$&$10$&$10$&$2$&$A_5, S_5, \PSL(2,9),S_6, \PGL(2,9), M_{10},\PGGL(2,9)$ (out of $9$)\\
8&$4$&$24$&$ 5\cdot 7\cdot 11$&$11$&$21$&$2$&$A_7,S_7$ (out of $9$)\\
10&$4$&$32$&$3^2\cdot 5\cdot 7$&$10$&$28$&$2$&$\PSU(3,3),\PGGU(3,3),\PSp(6,2)$ (out of $14$)\\
11&$4$&$32$&$3^2\cdot 5\cdot 7$&$28$&$10$&$4$&$\PGL(2,9),M_{10},\PGGL(2,9),A_{10},S_{10}$ (out of $9$)\\
\hline
\end{tabular}

\caption{Possibilities for the primitive group $D$
}\label{t:lastcases}
\end{center}
\end{table}

\begin{lemma}\label{le:Kneq1} Assume that Hypothesis \ref{h} holds, and that the parameters are as in one of the Lines 1, 7, 8, 10, 11 of Table~\ref{t:remainingcases}. Then 
\begin{enumerate}[(a)]
\item the kernel $K$ of the $G$-action on $\mathcal{C}$ is nontrivial;
\item     the group $D=G^\mathcal{C}$ is one of those listed in Table~\ref{t:lastcases};
\item either $K^\Delta$ is primitive, or  Line 11 holds with $L=\PGL(2,7)$ and $K^\Delta=\PSL(2,7)$;
\item $K$ acts faithfully on each  $\Delta\in \mathcal{C}$.

\end{enumerate}
\end{lemma}

\begin{proof} 

(a) In order to prove that $K\ne1$, we assume to the contrary that $K=1$. Then $G$ is isomorphic to $D$, a primitive group of degree $d$.
The fact that the number of flags $vr=bk$ must divide $|D|$ implies that $D=A_d$ or $S_d$ for Lines 1, 7, 8.  Applying Proposition \ref{le:Druleout} with the prime $p=7$ for Lines 1 and 7, or $p=13$ for Line  8 leads to contradictions.
We deal with Lines 10, 11 separately, considering the stabilisers $G_\Delta$ and $G_\alpha$, for $\alpha\in\Delta$.
%
%

In Line 11, $G\cong D$, a primitive subgroup of $S_{10}$ of order divisible by $bk=2^5.3^2.5.7$, and this implies that $G=A_{10}$ or $S_{10}$. Thus $G_\Delta$ is $A_{9}$ or $S_{9}$, respectively, and $G_\Delta$ has a subgroup $G_\alpha$ of index $c=28$, but neither $A_9$ nor $S_9$ has such a subgroup. 

Finally, in Line 10, $G\cong D$ is a primitive subgroup of $S_{28}$ of order divisible by $2^5.3^2.5.7$ but not by $p=17$ (by Proposition \ref{le:Druleout}). Hence, using \cite{MAGMA1} we see that $G=A_8$, $S_8$  or $\PSp(6,2)$, and so $G_\Delta=S_6$, $S_6\times 2$, or $\PSU(4,2):2$, respectively (in each case there is a unique conjugacy class of these subgroups of index $d=28$), and $G_\Delta$ has a subgroup $G_\alpha$ of index $c=10$. However in the third case there is no such subgroup, so  $G=A_8$ or $S_8$, and $G_\alpha = 3^2:D_8$ or $(3^2:D_8)\times 2$, respectively, (again unique up to conjugacy), and we identify the point set with the set of right cosets of $G_\alpha$ in $G$. Also a block stabiliser $G_B$ lies in a unique conjugacy class of subgroups of index $b=315$, and a computation with \cite{magma} shows that $G_B$ has $3$ or $2$ orbits of length $k=32$ on points, according as $D=A_8$ or $S_8$ respectively. However  a further computation with \cite{magma} (as described in Remark~\ref{rem:magma}) shows that, for each of these $32$-subsets $X$, the set of $G$-images of $X$ is not the block-set of a $2$-design. 
Thus in all cases we conclude that $K\neq 1$. 

(b) Since $L$ is primitive, the $K$-orbits on points have size $c$. We computationally exploit Proposition \ref{le:DBorbits}(b) with the appropriate  $x:=c/\gcd(c,\ell)$ to find the possibilities for $D$ listed in 
 Table~\ref{t:lastcases}; in each line of  Table~\ref{t:lastcases} we note the total number of primitive groups of degree $d$ (to indicate how many have been eliminated by these restrictions).

(c)
By Lemma \ref{le:prim}, $L$ is primitive. Since $K^\Delta\neq 1$, it follows that $K^\Delta$ is transitive and contains $\Soc(L)$. For $c=10$ or $c=11$, all primitive groups of degree $c$ have a primitive socle, so $K^\Delta$ is primitive in Lines 1, 7, 8, 10. In Line 11, we have
$c=28$, and there is a single  primitive group of degree $c$ with an  imprimitive socle, namely $\PGL(2,7)$. So  $K^\Delta$ is primitive also in Line 11, except if $L=\PGL(2,7)$ and $K^\Delta=\PSL(2,7)$.


(d)
Note that the group $G$ permutes the set $\{K_{(\Delta)}\mid \Delta\in\mathcal{C}\}$ by conjugation. 
For $\Delta'\ne\Delta$, since $K^{\Delta'}$ is transitive and $(K_{(\Delta)})^{\Delta'}\unlhd K^{\Delta'}$, it follows that either $(K_{(\Delta)})^{\Delta'}\ne 1$ with equal length orbits, or $K_{(\Delta)}$ fixes $\Delta'$ pointwise, and in the latter case $K_{(\Delta)}= K_{(\Delta')}$ (since these subgroups are $G$-conjugate). For a given $\Delta$, the subset $\{\Delta'\in\mathcal{C}\mid K_{(\Delta)}= K_{(\Delta')}\}$ is a block of imprimitivity for $G^\mathcal{C}=D$ containing $\Delta$, and since $D$ is primitive, it follows that either the subgroups $K_{(\Delta)}$ are all equal and hence are trivial $K_{(\Delta)}=1$, or the subgroups $K_{(\Delta)}$ are pairwise distinct and  the set of fixed points of $K_{(\Delta)}$ is precisely $\Delta$. Thus it is sufficient to prove that $K_{(\Delta)}=1$  for some $\Delta\in\mathcal{C}$.


Suppose, in order to derive a contradiction, that $K_{(\Delta)}\ne 1$. Then as explained above, for all $\Delta'\ne\Delta$,  the normal subgroup $K_{(\Delta)}^{\Delta'}$ of $K^{\Delta'}$ is nontrivial. Since $K^{\Delta'}$ is either primitive or a nonabelian simple group, by part (c), this implies that $K_{(\Delta)}^{\Delta'}$ is transitive. This leads to a contradiction as follows: each pair $\{\alpha,\beta\}\subseteq \Delta$ is fixed by $K_{(\Delta)}$, and hence $K_{(\Delta)}$ leaves invariant the set of $\lambda$ blocks containing $\{\alpha,\beta\}$. Let $\Delta'\neq \Delta$ be a class intersecting non-trivially at least one of these blocks. Then $K_{(\Delta)}$ leaves invariant the (non-empty) union of the subsets $B\cap \Delta'$ over the $\lambda$ blocks $B$ containing $\{\alpha,\beta\}$. This non-empty union is a subset of size at most $\lambda \ell<c=|\Delta'|$, and this contradicts the fact that $K_{(\Delta)}^{\Delta'}$ is transitive. Hence  $K_{(\Delta)}= 1$, and part (d) is proved.
\end{proof}

\medskip
We now give ad hoc arguments, involving many of the results from Section~\ref{sec:methods}, to deal with the five remaining Lines of  Table~\ref{t:remainingcases}, namely the Lines of Table~\ref{t:lastcases}. 

\begin{lemma}\label{cor:line9}
There are no examples with parameter tuple in Line 8  of Table~\ref{t:remainingcases}.
\end{lemma}
\begin{proof}
In Line 8 we have $(c,d,\ell)=(11,21,2)$ and $\lambda=4$. Note that $p=5>\lambda$ and $5$ does not divide $c$. If $5$ divides $|K|$, then by Proposition \ref{cor:subdesign}, there exists an integer $f$ such that $2< f<11$, $f\equiv 11\equiv 1\pmod 5$, and $k-1=23$ divides $\lambda(df-1)$, which is a contradiction. Hence $5$ does not divide $|K|$.
Now by Lemma \ref{le:Kneq1}, $K\cong K^\Delta$ is a non-trivial normal subgroup of $L$, so $K^\Delta$ contains the socle of $L$. Moreover, by Proposition \ref{le:L}(a0, $L$ is 2-transitive of degree $c=11$. The socles of such groups are $C_{11}$, $\PSL(2,11)$, $M_{11}$, or $A_{11}$, and the only one with order coprime to $5$ is $C_{11}$. Thus  $C_{11}\unlhd K^\Delta\unlhd L=\AGL(1,11)$. Since  $5$ does not divide $|K|$, we have $K\cong K^\Delta=C_{11}$ or $C_{11}\rtimes C_2$. 

Let  $K_0:=O_{11}(K)\cong C_{11}$, and let  $C:=C_G(K_0)$. By Lemma~\ref{le:Kneq1}(b), $D=G/K$ is $A_7$ or $S_7$. Since the conjugation action of $G$ on $K_0$ induces a subgroup of $\Aut(K_0)\cong C_{10}$, it follows that $C^\mathcal{C}=CK/K$ contains $A_7$. Thus $|G:CK|\leq 2$, and as $C\cap K=K_0$, we also have $|CK:C|=|K:C\cap K|\leq 2$. Hence $|G:C|$ divides $4$. However this is a contradiction since, on the one hand, $G/C$ is isomorphic to the subgroup of $\Aut(K_0)$ induced by $G$, and on the other hand, $G_\Delta^\Delta=\AGL(1,11)$ induces the whole group $\Aut(K_0)\cong C_{10}$ on $K_0$.
\end{proof}

\begin{lemma}\label{cor:lines1and8and11}
There are no examples with parameter tuple in Line $1, 7$ or $10$  of Table~\ref{t:remainingcases}.
\end{lemma}
\begin{proof}
In these Lines, $(c, \ell)=(10, 2)$, and the triple $(\lambda, d,k)$ is $(3, 10,12), (4, 10,12)$ or $(4,28,32)$ in  Line 1, 7 or 10, respectively. By Lemma~\ref{le:Kneq1}, for a class $\Delta\in\mathcal{C}$,  $K^\Delta \cong K$  is primitive of degree $c=10$, and the possibilities for $\Soc(K)$ are $A_5, A_6$ and $A_{10}$. Also $\Soc(L)\unlhd K^\Delta\unlhd L$, and $L$ is $2$-transitive by Proposition~\ref{le:L}(a), so only $A_6$ and $A_{10}$ are possible for $\Soc(K)$. In particular $\Soc(K^\Delta)$ is $2$-transitive.  Let $P$ be a Sylow $3$-subgroup of $K$. Since $K^\Delta$ is $2$-transitive, it follows that $P^\Delta$ has orbits of lengths $1,9$. It follows that $P$ fixes a unique point in each of the classes of $\mathcal{C}$. Let $\alpha, \beta$ be distinct points fixed by $P$ with $\alpha\in\Delta$. Suppose first that $\lambda=4$. Then by Proposition~\ref{le:plam}, since $\lambda\equiv 1\pmod{3}$, $P$ fixes setwise at least one block containing $\{\alpha,\beta\}$, say $P$ fixes $B$ setwise. This implies that $P$ fixes setwise the $\ell$-subset $B\cap \Delta$, which is a contradiction since $\ell=2$ and $P$ fixes a unique point of $\Delta$. Thus we are in Line 1 with $\lambda=3$, and from the argument just given we may assume that $P$ fixes none of the three blocks $B_1, B_2, B_3$ containing $\{\alpha,\beta\}$. Thus $P$ permutes these three blocks cyclically. Noting that $\ell=2$, let $B_i\cap \Delta=\{\alpha, \alpha_i\}$. Then $P$ fixes $\{\alpha_1,\alpha_2,\alpha_3\}$ setwise, contradicting the fact that $P^\Delta$ has orbits of lengths $1,9$. 
This contradiction completes the proof.
\end{proof}

\begin{lemma}\label{cor:line12}
There are no examples with parameter tuple in Line 11 of Table~\ref{t:remainingcases}.
\end{lemma}

\begin{proof}
Here $(c, d, \lambda,\ell, k) = (28, 10, 4, 4, 32)$ and the number of blocks is $b=3^2\cdot 5\cdot 7$. By Lemma \ref{le:prim}, the groups $L, D$ are primitive of degree $c, d$ respectively.
By Proposition~\ref{le:L}(b), for distinct $\alpha\beta\in\Delta$,  the setwise stabiliser 
$L_{\{\alpha,\beta\}}$ has an orbit in $\Delta\setminus\{\alpha,\beta\}$ of size at most $(\ell-2)\lambda=8$. Checking with \cite{MAGMA1} we find that, of the $14$ primitive groups $L$ of degree $28$,  this property holds for only $7$ of them. Thus $L$ is one of 
\[
\PGL(2,7), \PSL(2,8),\PGGL(2,8),\PSU(3,3), \PGGU(3,3), A_8, S_8.
\]
In all cases $\Soc(L)$ is a nonabelian simple group, and is not regular, and the centraliser $C_{{\rm Sym}(\Delta)}(\Soc(L))=1$.
Also, by Lemma~\ref{le:Kneq1},  $K^\Delta \cong K\ne 1$ and hence $\Soc(L)\leq K^\Delta\unlhd L$.

Let $C:=C_G(K)$. Then $C\cap K\leq \prod_{\Delta\in\mathcal{C}}C_{{\rm Sym}(\Delta)}(K^\Delta)$, and it follows that $C\cap K=1$. Now the conjugation action of $G$ on $K$ induces a subgroup of $\Aut(K)$ isomorphic to $G/C$, and the subgroup $CK/C \cong K$ induces the group of inner automorphisms  of $K$. Hence $G/CK$ is isomorphic to a section of the outer automorphism group of $\Soc(L)$, which in all cases is a group of order at most $3$. Thus  $C^\mathcal{C}=CK/K$ has index at most $3$ in $D=G^\mathcal{C}=G/K$, and we note that  $C^\mathcal{C}=CK/K\cong C$. It follows that $\Soc(D)\leq C^\mathcal{C}\unlhd D$, and by Lemma~\ref{le:Kneq1}(b),  $\Soc(D)\in\{ A_6, A_{10}\}$. 

We claim that $C$ has $c=28$ orbits of length $d=10$  in $\mathcal{P}$. Since $C\lhd G$, $C$ has equal length orbits in $\mathcal{P}$, and since $C^\mathcal{C}$ contains $\Soc(D)$, $C^\mathcal{C}$ is transitive. Thus the  length of the $C$-orbits is $10\cdot c_0$, where $c_0$ is the (constant) length of the $C_\Delta$-orbits in $\Delta$. Now $C_\Delta^\Delta$ is a  normal subgroup  of the primitive group $G_\Delta^\Delta$ of degree $c=28$. If $C_\Delta^\Delta\ne 1$ then $C_\Delta^\Delta$ is transitive and in particular $7$ divides $|C_\Delta|$, which divides $|D|$. Thus $\Soc(D)=A_{10}$, and hence $C_\Delta=A_9$ or $S_9$.  Checking with Magma \cite{magma} we find that $C_\Delta$ has no subgroup of index $28$. Thus $C_\Delta^\Delta= 1$, and this implies that the $C$-orbits have length $10\cdot c_0=10$, proving the claim. 
The $C$-orbits form a $G$-invariant partition of $\mathcal{P}$ consisting of $28$ classes of size $10$. Therefore the design (if it existed) would also arise as a design with the parameters of Line 10. However there are no examples for Line 10 by Lemma~\ref{cor:lines1and8and11}.   
\end{proof}

\medskip
The collection of lemmas in this section together prove Theorem~\ref{t:lam34}.

\end{document}